\documentclass[12pt, leqno, draft]{article}

\input epsf
\epsfverbosetrue

\begin{document}

\title{What is a metric space?}

\author{Stephen Semmes	\\
	Rice University}

\date{}

\maketitle

\begin{abstract}
A \emph{metric space} is a mathematical notion that includes classical
Euclidean geometry and a variety of other situations.  Some basic
examples and their properties are briefly discussed.
\end{abstract}

\tableofcontents

\section{The Euclidean plane}
\label{euclidean plane}

\begin{figure}

\makebox[\textwidth]{
\epsfxsize=200pt
\epsfbox{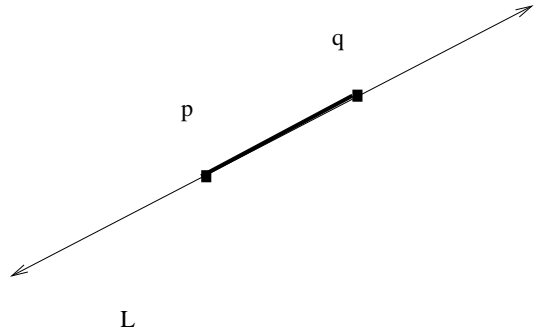}
}

\caption{\label{line} Two points $p$, $q$ on a line $L$, and the line
segment that connects them.}

\end{figure}

	Let $p$ and $q$ be a pair of distinct points in a Euclidean
plane.  There is a unique line $L$ in the plane passing through $p$
and $q$.  The \emph{distance} $d(p, q)$ between $p$ and $q$ may be
defined as the length of the line segment on $L$ that connects $p$ to
$q$, as in Figure \ref{line}.

\begin{figure}

\makebox[\textwidth]{
\epsfxsize=200pt
\epsfbox{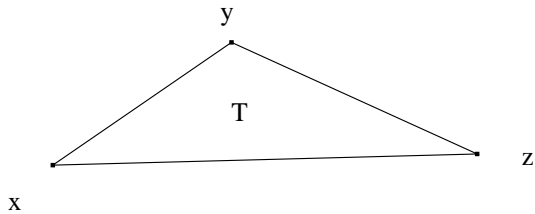}
}

\caption{\label{triangle} A triangle $T$ with vertices $x$, $y$, and $z$.}

\end{figure}

	If $x$, $y$, and $z$ are any three points in the plane, then
\begin{equation}
\label{triangle inequality}
	d(x, z) \le d(x, y) + d(y, z).
\end{equation}
This is known as the \emph{triangle inequality}.  If $T$ is the
triangle with vertices $x$, $y$, $z$, as in Figure \ref{triangle},
then (\ref{triangle inequality}) says that the length of the side of
$T$ connecting $x$ to $z$ is less than or equal to the sum of the
lengths of the other two sides.

\section{Abstract metric spaces}
\label{abstract metric spaces}

	A \emph{metric space} is a set $M$, which is to say a
collection of some sort of objects, for which the distance $d(p, q)$
between any two elements $p$ and $q$ of $M$ has been defined and
satisfies certain conditions.  Specifically,
\begin{equation}
	d(q, p) = d(p, q) \ge 0
\end{equation}
for every $p, q \in M$, with $d(p, q) = 0$ exactly when $p = q$, and
the triangle inequality (\ref{triangle inequality}) should hold for
every $x, y, z \in M$.  Thus, ordinary Euclidean geometry on a plane
is an example of a metric space.  The distance function $d(p, q)$ is
also known as a \emph{metric}.  On any set $M$, a metric can be
defined by putting $d(p, q) = 1$ when $p \ne q$, called the
\emph{discrete metric}.  If $(M, d(x, y))$ is a metric space, and $E$
is a subset of $M$, then the restriction of $d(x, y)$ to $x, y \in E$
defines a metric on $E$.  The triangle inequality and other
requirements of a metric are automatically satisfied on $E$, since
they hold on $M$.

\section{The real line}
\label{section on real line}

\begin{figure}

\makebox[\textwidth]{
\epsfxsize=250pt
\epsfbox{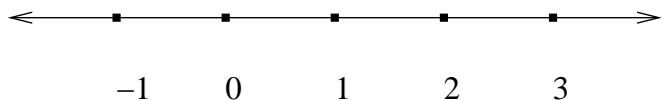}
}

\caption{\label{real line} The real line.}

\end{figure}

	The \emph{real line} ${\bf R}$ consists of all real numbers,
as in Figure \ref{real line}.  If $r$ is a real number, then the
\emph{absolute value} of $r$ is denoted $|r|$ and defined by $|r| = r$
when $r \ge 0$ and $|r| = - r$ when $r \le 0$.  For every $r, t \in
{\bf R}$,
\begin{equation}
\label{|r + t| le |r| + |t|}
	|r + t| \le |r| + |t|.
\end{equation}
More precisely, $|r + t| = |r| + |t|$ when $r$ or $t$ is equal to $0$,
or when $r$ and $t$ have the same sign, and otherwise there is some
cancellation in $r + t$.  The standard distance between $r$ and $t$
can be defined by $d(r, t) = |r - t|$.  In this case, the triangle
inequality reduces to (\ref{|r + t| le |r| + |t|}).  For if $x$, $y$,
and $z$ are arbitrary real numbers, then $x - z = (x - y) + (y - z)$,
and therefore
\begin{equation}
	|x - z| \le |x - y| + |y - z|.
\end{equation}
With this metric, the real line is equivalent to a line in Euclidean
geometry.

\section{Back to the plane}
\label{back to plane}

\begin{figure}

\makebox[\textwidth]{
\epsfxsize=250pt
\epsfbox{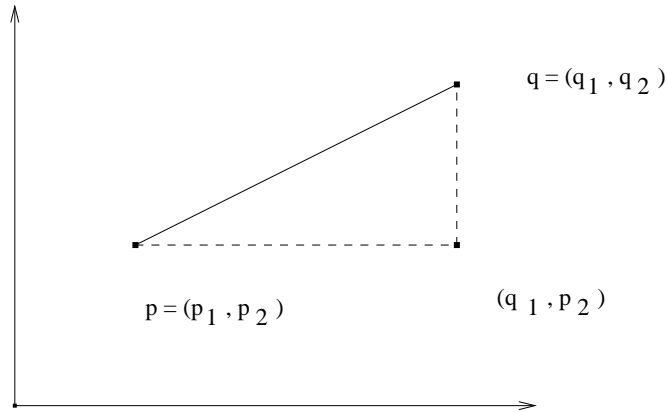}
}

\caption{\label{cartesian} A right triangle with vertices $p =
(p_1, p_2)$, $q = (q_1, q_2)$, and $(q_1, p_2)$.}

\end{figure}

	Points in a plane can also be described in terms of numbers,
using \emph{Cartesian coordinates}.  More precisely, a point $p$ in
the plane can be represented by two real numbers $p_1$ and $p_2$, the
first and second coordinates of $p$.  Suppose that
\begin{equation}
\label{p = (p_1, p_2) and q = (q_1, q_2)}
	p = (p_1, p_2) \quad\hbox{and}\quad q = (q_1, q_2)
\end{equation}
are two points in the plane represented by pairs of numbers in this
manner.  The usual Euclidean distance $d(p, q)$ between $p$ and $q$ is
determined by the formula
\begin{equation}
	d(p, q)^2 = (p_1 - q_1)^2 + (p_2 - q_2)^2.
\end{equation}
This follows from the well-known \emph{Pythagorean theorem}, which
states that the length of the hypotenuse of a right triangle is equal
to the sum of the squares of the lengths of the other two sides.
Remember that the hypotenuse of a right triangle is the side opposite
from the right angle.  Here we can consider the triangle in the plane
whose vertices are $p$, $q$, and the point corresponding to $(q_1,
p_2)$, as in Figure \ref{cartesian}.  This is a right triangle whose
hypotenuse is the line segment that connects $p$ to $q$, and whose
other two sides are parallel to the coordinate axes and have lengths
$|p_1 - q_1|$ and $|p_2 - q_2|$.

	One can also look at the distance between $p$ and $q$ as in
(\ref{p = (p_1, p_2) and q = (q_1, q_2)}) defined by
\begin{equation}
\label{d'(p, q) = |p_1 - q_1| + |p_2 - q_2|}
	d'(p, q) = |p_1 - q_1| + |p_2 - q_2|.
\end{equation}
This satisfies the requirements of a metric space too, just like the
Euclidean distance.  This is not the same as the Euclidean distance,
because the distance between $(0, 0)$ and $(1, 1)$ is $\sqrt{2}$ in
the Euclidean metric, and is equal to $2$ with respect to (\ref{d'(p,
q) = |p_1 - q_1| + |p_2 - q_2|}).  Another example is given by
\begin{equation}
\label{d''(p, q) = max (|p_1 - q_1|, |p_2 - q_2|)}
	d''(p, q) = \max (|p_1 - q_1|, |p_2 - q_2|),
\end{equation}
which is to say that
\begin{equation}
	d''(p, q) = |p_1 - q_1| \quad\hbox{when } |p_1 - q_1| \ge |p_2 - q_2|,
\end{equation}
and
\begin{equation}
	d''(p, q) = |p_2 - q_2| \quad\hbox{when } |p_1 - q_1| \le |p_2 - q_2|.
\end{equation}
One can check that this satisfies the requirements of a metric space
as well.  Note that the distance between $(0, 0)$ and $(1, 1)$ is
equal to $1$ with respect to (\ref{d''(p, q) = max (|p_1 - q_1|, |p_2
- q_2|)}).  The distance between $(0, 0)$ and $(1, 0)$ is equal to $1$
with respect to each of these three metrics.

\begin{figure}

\makebox[\textwidth]{
\epsfxsize=100pt
\epsfbox{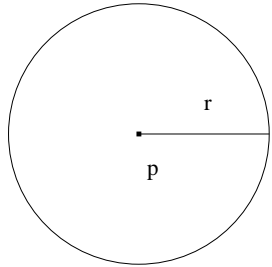}
}

\caption{\label{circle} An ordinary Euclidean circle with center $p$
and radius $r$.}

\end{figure}

\begin{figure}

\makebox[\textwidth]{
\epsfxsize=100pt
\epsfbox{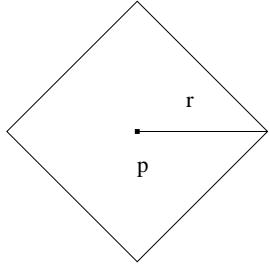}
}

\caption{\label{diamond} The set of points $q$ such that $d'(p, q) = r$
has the shape of a diamond.}

\end{figure}

\begin{figure}

\makebox[\textwidth]{
\epsfxsize=100pt
\epsfbox{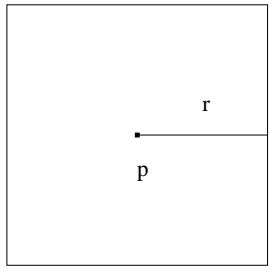}
}

\caption{\label{square} The set of points $q$ such that $d''(p, q) = r$
has the shape of a square.}

\end{figure}

	Let $p$ be a point in the plane, and let $r$ be a positive
real number.  The usual Euclidean circle in the plane with center $p$
and radius $r$ consists of the points $q$ in the plane such that
\begin{equation}
	d(p, q) = r,
\end{equation}
where $d(p, q)$ is the standard Euclidean distance, as in Figure
\ref{circle}.  By contrast, the set of points $q$ in the plane such
that
\begin{equation}
	d'(p, q) = r
\end{equation}
has the shape of a diamond, as in Figure \ref{diamond}.  Similarly,
the set of $q$ such that
\begin{equation}
	d''(p, q) = r
\end{equation}
is a square, as in Figure \ref{square}.

\begin{figure}

\makebox[\textwidth]{
\epsfxsize=100pt
\epsfbox{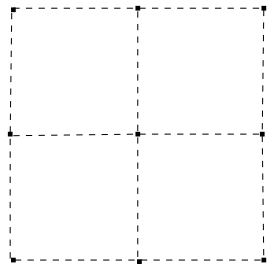}
}

\caption{\label{taxicab} The taxicab metric measures distance between
points in a grid by following horizontal and vertical segments.}

\end{figure}

        The metric $d'(p, q)$ is sometimes called the ``taxicab
metric'', because it measures the minimum distance from $p$ to $q$
along a combination of horizontal and vertical segments, with no
diagonal shortcuts.  More precisely, it is customary to consider
points whose coordinates are whole numbers, and to think of horizontal
and vertical lines through these points as being like streets, as in
Figure \ref{taxicab}.  Note that there can be more than one such path
of minimal length between two points.

\begin{figure}

\makebox[\textwidth]{
\epsfxsize=200pt
\epsfbox{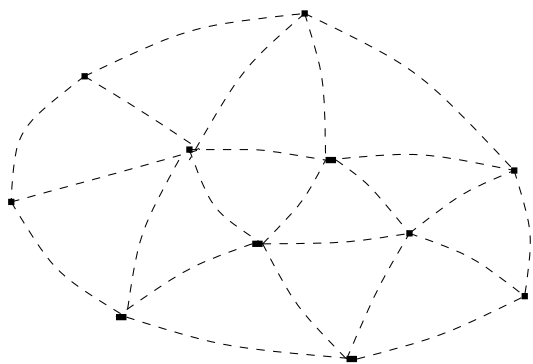}
}

\caption{\label{grid} A collection of points connected by various paths.}

\end{figure}

        The same idea can be applied to any collection of points
connected by various paths, as in Figure \ref{grid}.  If $x$, $y$, and
$z$ are three of the points, then any route from $x$ to $y$ can be
combined with a route from $y$ to $z$ to get a route from $x$ to $z$.
This means that the distance from $x$ to $z$ is less than or equal to
the sum of the distances from $x$ to $y$ and from $y$ to $z$, since
this distance is defined by minimizing the lengths of paths between
two points.  Thus the triangle inequality holds automatically in this
type of situation.

\section{Spheres and other surfaces}
\label{spheres and surfaces}

        Suppose now that $M$ is a nice round two-dimensional sphere in
ordinary three-dimensional Euclidean space.  One way to measure
distance between points on $M$ is to simply use the Euclidean distance
from the ambient three-dimensional space.  Another way is to minimize
lengths of paths on $M$.  If $p$ and $q$ are distinct points on $M$,
then the Euclidean distance between $p$ and $q$ is the same as the
length of the line segment connecting $p$ and $q$, which is not
contained in $M$.  It turns out that the path on $M$ connecting $p$
and $q$ with minimal length is the short arc of a \emph{great circle}
passing through $p$ and $q$.  A great circle on a sphere is a circle
with the same radius as the sphere, like longitudinal circles through
the poles, or the equator.  Equivalently, each circle on a sphere is
the intersection of the sphere with a plane, and a great circle is the
intersection of the sphere with a plane passing through the center of
the sphere.  If $p$ and $q$ are antipodal points on $M$, then there
are infinitely many great circles passing through $p$ and $q$.  The
arcs of these great circles connecting $p$ and $q$ have the same
length, which is one-half the circumference of the sphere.  Otherwise,
there is exactly one great circle on $M$ passing through $p$ and $q$,
and the shorter arc of this great circle connecting $p$ and $q$ is the
path of minimal length on $M$ connecting $p$ and $q$.

        More generally, if $M$ is a nice surface in three-dimensional
Euclidean space, then one can consider the distance defined by
minimizing lengths of paths on $M$ between a given pair of points.
Instead of a flat plane or a round sphere, $M$ could be bumpy, like a
golf ball.  The shape of a path of minimal length depends on the
geometry of the surface, and may not be so simple as a line segment,
or a circular arc, or the intersection of $M$ with a plane.  There may
be many paths of minimal length, as for antipodal points on a sphere.

\begin{figure}

\makebox[\textwidth]{
\epsfxsize=100pt
\epsfbox{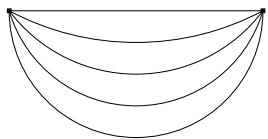}
}

\caption{\label{arcs} Circular arcs and the line segment connecting a
pair of points in a Euclidean plane.}

\end{figure}

        Implicit in the discussion of the case of a sphere is a
geometric fact that makes sense already in two dimensions: the length
of the short arc of a circle between a fixed pair of points in a
Euclidean plane decreases as the radius of the circle increases, as in
Figure \ref{arcs}.  One might also say that the curvature of the
circle decreases as the radius increases, and that a line is like a
circle with infinite radius and curvature zero.  The same principle
applies to circular arcs in three dimensions, by rotating an arc so
that it lies in the same plane as another arc.

        Why exactly does the length of a circular arc between two
fixed points decrease as the radius of the circle increases?  Here are
a few ways to look at this question.  It does not matter where the
endpoints of the arcs are, and so we can suppose that they correspond
to $-1$ and $1$ on the $x$-axis.  Each circular arc can be described
as the graph of a real-valued function on the interval $[-1, 1]$, and
one can check that the magnitude of the slope of the tangent to a
point on such an arc over a particular point in the interval decreases
as the radius of the circle increases.  This is pretty clear from the
picture, and it implies that the length of the arc decreases as the
radius increases.

        Alternatively, we can reformulate the original statement by
saying that the distance between two points on the unit circle divided
by the length of the shorter arc between them increases as the points
move closer together.  One can see that this is equivalent to the
previous version by rescaling the picture.  The length of the shorter
arc determines an angle $\theta$, $0 < \theta \le \pi$, and the
distance between the two points is equal to $2 \sin (\theta/2)$.  Thus
one would like to show that
\begin{equation}
        \frac{\sin t}{t}
\end{equation}
is monotone decreasing when $0 < t \le \pi/2$, which can be treated as
an exercise in calculus.  Note that $t = 0$ can be included using the
continuous extension of the ratio equal to $1$ there.

        As a another approach, let $C$ be a circle in the plane with
center $z$, and let $E$ be the set of points on $C$ or in the exterior
of $C$.  By definition, the \emph{circular projection} sends a point
$p \ne z$ in the plane to the element of $C$ on the line segment
connecting $p$ to $z$.  If $p$ and $q$ are elements of $E$, then the
distance between their projections in $C$ is less than or equal to
their distance in the plane.  In terms of vector calculus, the
differential of the circular projection at any point $p \in E$
corresponds to a linear mapping that does not increase the Euclidean
norm of a vector.  This implies that the length of any curve in $E$
with endpoints on $C$ is greater than or equal to the length of the
shorter arc of $C$ with the same endpoints, and applies in particular
to circular arcs with endpoints on $C$ associated to circles of
smaller radius.

        Let $M$ be the unit sphere in three-dimensional Euclidean
space, consisting of points whose distance to the origin is equal to
$1$, and let $d(p, q)$ be the ordinary Euclidean distance between any
two points $p, q \in M$.  For each $p, q \in M$, $0 \le d(p, q) \le
2$, and there is a unique $\widetilde{d}(p, q) \in [0, \pi]$ such that
\begin{equation}
        \sin \Big(\frac{\widetilde{d}(p, q)}{2}\Big) = \frac{d(p, q)}{2}.
\end{equation}
Equivalently, $\widetilde{d}(p, q)$ is the length of the shorter arc
of a great circle on the unit sphere passing through $p$ and $q$.  To
show that this defines a metric on $M$, we need to check that the
triangle inequality
\begin{equation}
        \widetilde{d}(p, r) \le \widetilde{d}(p, q) + \widetilde{d}(q, r)
\end{equation}
holds for every $p, q, r \in M$.  Consider the set $A$ of $r' \in M$
with
\begin{equation}
        \widetilde{d}(q, r') = \widetilde{d}(q, r),
\end{equation}
which is the same as
\begin{equation}
        d(q, r') = d(q, r).
\end{equation}
Thus $A$ is a circle, which reduces to a single point in the trivial
cases where $q$ and $r$ are equal or antipodal.  Let $r_0$ be the
point in $A$ such that
\begin{equation}
        d(p, r') \le d(p, r_0)
\end{equation}
for each $r' \in A$, which is the same as
\begin{equation}
        \widetilde{d}(p, r') \le \widetilde{d}(p, r_0)
\end{equation}
for every $r \in A$.  One can show that $p$, $q$, and $r_0$ lie on a
great circle, as in the next paragraph.  The triangle inequality holds
for points on a great circle, which is to say that
\begin{equation}
        \widetilde{d}(p, r_0) \le \widetilde{d}(p, q) + \widetilde{d}(q, r_0).
\end{equation}
Therefore
\begin{equation}
        \widetilde{d}(p, r) \le \widetilde{d}(p, r_0)
         \le \widetilde{d}(p, q) + \widetilde{d}(q, r_0)
           = \widetilde{d}(p, q) + \widetilde{d}(q, r),
\end{equation}
as desired.  One can also check that equality holds in the triangle
inequality only for points on a great circle.

        In a Euclidean plane, the distance between a point $x$ and
elements of a circle is maximized and minimized on the intersection of
the circle with the line through $x$ and the center of the circle.
The same holds for a point and a sphere in three-dimensional Euclidean
space.  Suppose that we have a point $x$ and a circle in
three-dimensional Euclidean space, where $x$ is not on the plane of
the circle.  Let $\widetilde{x}$ be the projection of $x$ onto the
plane of the circle, which is to say that the line through $x$ and
$\widetilde{x}$ is perpendicular to the plane of the circle.  The
distance from $x$ to an element of the circle can be expressed in
terms of the distance from $x$ to $\widetilde{x}$ and the distance
from $\widetilde{x}$ to the point in question, by the Pythagorean
theorem.  Maximizing or minimizing the distance from the circle to $x$
is hence the same as maximizing or minimizing the distance to
$\widetilde{x}$, respectively.  The maximum and minimum occur on the
intersection of the circle with the line through $\widetilde{x}$ and
the center of the circle, since $\widetilde{x}$ is on the plane of the
circle.  In terms of $x$, the maximum and minimum occur on the plane
passing through $x$ and the center of the circle which is
perpendicular to the circle.  For the problem in spherical geometry in
the previous paragraph, the Euclidean center of the circle $A$ is on
the line passing through $q$ and the origin, and the plane of $A$ is
perpendicular to this line.  The plane that contains the origin, $p$,
and $q$ also contains the Euclidean center of $A$ and is perpendicular
to the plane of $A$, and so the distance from $p$ to elements of $A$
is maximized and minimized on the intersection of this plane with $A$.
Thus the maximum and minimum occur on a great circle through $p$ and
$q$.

        By construction,
\begin{equation}
         d(p, q) \le \widetilde{d}(p, q)
\end{equation}
for every $p, q \in M$.  The distances are practically the same
locally, in the sense that for every $\epsilon > 0$ there is a $\delta
> 0$ such that
\begin{equation}
         \widetilde{d}(p, q) \le (1 + \epsilon) \, d(p,q)
\end{equation}
when $d(p, q) < \delta$, because
\begin{equation}
         \lim_{x \to 0} \sin x / x = 1.
\end{equation}
This means that the two metrics are the same ``infinitesimally'' at
each point, and in particular that the length of a path in $M$ is the
same relative to either metric.  One can use this to show that arcs of
great circles with length $\le \pi$ are the paths of minimal length,
since the length of a path is always greater than or equal to the
distance between its endpoints.  Hence the spherical metric
$\widetilde{d}(p, q)$ is the same as the metric defined by minimizing
lengths of paths in this case.

\section{Intrinsic and extrinsic geometry}
\label{intrinsic, extrinsic geometry}

\begin{figure}

\makebox[\textwidth]{
\epsfxsize=250pt
\epsfbox{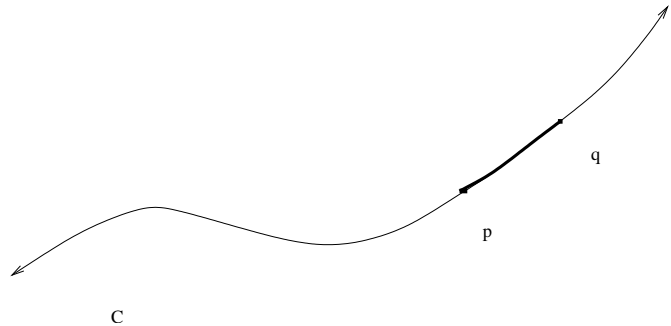}
}

\caption{\label{curve} A simple curve $C$ in the plane, and the arc
connecting two points $p$ and $q$ on $C$.}

\end{figure}

        Let $C$ be a nice smooth curve in the plane, as in Figure
\ref{curve}.  Suppose that $C$ is \emph{simple}, i.e., without
crossings.  It may be that $C$ is bounded, like a line segment or
circular arc, or unbounded, like a line or ray.  As in the previous
section, the distance between two points $p$ and $q$ on $C$ might be
measured using the Euclidean metric from the plane, or using the
length of the arc on $C$ between $p$ and $q$.  This is the path of
minimal length on $C$ connecting $p$ to $q$, which is more complicated
on a surface, since the paths can move around more.  As usual, the
Euclidean distance between $p$ and $q$ is less than or equal to the
length of the arc on $C$ between them.  Normally the inequality is
strict, because of curvature.  With respect to the distance defined by
the length of arcs, $C$ looks flat.  Using a parameterization of $C$
by arc length, as in vector calculus, $C$ looks the same as if it were
contained in a line relative to this distance.  Thus $C$ may be
intrinsically flat even if it is extrinsically curved.  However,
surfaces are typically curved intrinsically as well as extrinsically.

\section{Non-Euclidean geometry}
\label{non-euclidean}

        If $L$ is a line in a Euclidean plane and $p$ is a point in
the plane not on $L$, then there is exactly one line $L'$ in the plane
passing through $p$ and parallel to $L$.  This is the famous
\emph{parallel postulate} in Euclidean geometry.  A well-known
question asked whether the parallel postulate could be derived from
the other postulates.  Eventually it was discovered that this is not
possible, because there are non-Euclidean geometries in which the
parallel postulate does not work while the other postulates hold.  One
of these is based on spherical geometry, with great circles playing
the role of lines.  More precisely, it is better to use projective
space, where antipodal points in the sphere are identified.  In this
case, the problem with the parallel postulate is that great circles
always intersect.  In \emph{hyperbolic geometry}, the problem with the
parallel postulate is that there is more than one line through a point
$p$ that is parallel to a fixed line not containing $p$.  Hyperbolic
geometry can also be represented by a metric space, where arcs of
``lines'' minimize length and determine distance.

\section{Symmetry}
\label{symmetry}

         A basic question that one might ask about a metric space is
what kind of symmetry it has.  On the real line, for example, there
are symmetries by translation.  For each real number $a$, one can
translate every element of ${\bf R}$ by $a$ without changing the
distances.  For if $\tau_a(x) = x + a$, then
\begin{equation}
         |\tau_a(x) - \tau_a(y)| = |(x + a) - (y + a)| = |x - y|
\end{equation}
for every $x, y \in {\bf R}$.  Similarly, if $\rho(x) = -x$, then
\begin{equation}
         |\rho(x) - \rho(y)| = |(-x) - (-y)| = |y - x| = |x - y|
\end{equation}
for every $x, y \in {\bf R}$.  This is the reflection of ${\bf R}$
about $0$, and the reflection of ${\bf R}$ about any other point
preserves distances too.  The reflection about any other point can
also be expressed in terms of $\rho$ and suitable translations.

         Translation symmetries on the plane can be described similarly
in coordinates by
\begin{equation}
         (x_1, x_2) \mapsto (x_1 + a_1, x_2 + a_2),
\end{equation}
where $a_1, a_2 \in {\bf R}$.  Reflection about the origin can be
given by
\begin{equation}
         (x_1, x_2) \mapsto (-x_1, -x_2).
\end{equation}
One can reflect about the axes independently, as in
\begin{equation}
         (x_1, x_2) \mapsto (-x_1, x_2)
\end{equation}
and
\begin{equation}
         (x_1, x_2) \mapsto (x_1, - x_2),
\end{equation}
or about the line $x_1 = x_2$, as in
\begin{equation}
         (x_1, x_2) \mapsto (x_2, x_1).
\end{equation}
A remarkable feature of Euclidean geometry is that the metric is
preserved by arbitrary rotations.  This does not work for the metrics
$d'(p, q)$ and $d''(p, q)$ defined in Section \ref{back to plane},
although translations and the reflections just mentioned preserve
these metrics.

         For a two-dimensional round sphere in three-dimensional
Euclidean space, there are symmetries by rotation about the center of
the sphere.  One can reflect about the center of the sphere, or about
any plane through the center of the sphere.  As in the previous
examples, these symmetries are sufficient to move any point in the
space to any other point in the space.

\section{A nice little picture}
\label{nice little picture}

        Let $(M, d(x, y))$ be a metric space.  For each $p \in M$ and
$r > 0$, $B(p, r)$ denotes the \emph{open ball} in $M$ with center $p$
and radius $r$, defined by
\begin{equation}
        B(p, r) = \{x \in M : d(p, x) < r\}.
\end{equation}
In the real line with the standard metric, this is the same as the
open interval $(p - r, p + r)$.  In a Euclidean plane, $B(p, r)$ is a
round disk.

\begin{figure}

\makebox[\textwidth]{
\epsfxsize=200pt
\epsfbox{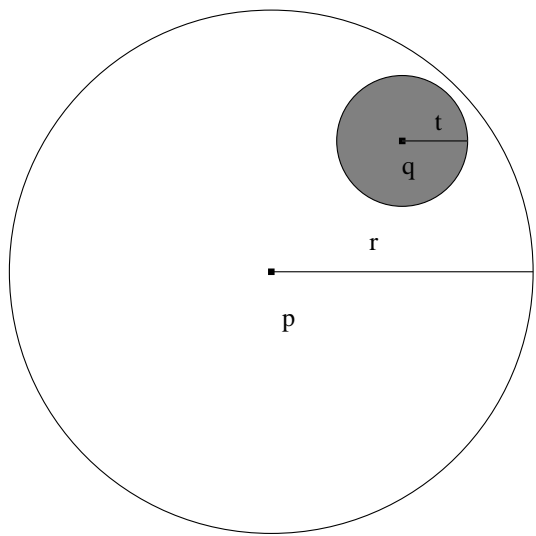}
}

\caption{\label{inclusion} If $d(p, q) < r$ and $0 < t \le r - d(p,
q)$, then $B(q, t) \subseteq B(p, r)$.}

\end{figure}

        If $q \in B(p, r)$ and $0 < t \le r - d(p, q)$, then
\begin{equation}
        B(q, t) \subseteq B(p, r),
\end{equation}
as in Figure \ref{inclusion}.  For if $x$ is any element of $B(q, t)$,
then $d(q, x) < t$, and
\begin{equation}
        d(p, x) \le d(p, q) + d(q, x) < d(p, q) + t \le r,
\end{equation}
by the triangle inequality.  Thus $d(p, x) < r$, as desired.  More
precisely, Figure \ref{inclusion} shows a Euclidean plane and is
suggestive of the general case.

\end{document}